\documentclass[12pt]{article}

\usepackage{pstricks,pst-node}
\usepackage{amsmath,amsthm}
\usepackage{amssymb}
\usepackage{color}
\usepackage{xspace}
\usepackage[colorlinks=true,
linkcolor=green,
filecolor=brown,
citecolor=green]{hyperref}

\begin{document}
\newtheorem{theorem}{Theorem}
\newtheorem{defn}[theorem]{Definition}
\newtheorem{lemma}[theorem]{Lemma}
\newtheorem{prop}[theorem]{Proposition}
\newtheorem{cor}[theorem]{Corollary}
\begin{center}
{\Large
Flexagons Yield a Curious Catalan Number Identity                           \\ 
}

\vspace{10mm}
DAVID CALLAN  \\
Department of Statistics  \\
\vspace*{-1mm}
University of Wisconsin-Madison  \\
{\bf callan@stat.wisc.edu}  \\
\vspace*{1mm}
June 2, 2010 \\
\vspace*{5mm}
\small{\textbf{Dedicated to the memory of Martin Gardner, 1914--2010.}}
\end{center}

\begin{abstract}
Hexaflexagons were popularized by the late Martin Gardner in his first 
Scientific American column in 1956. 
Oakley and Wisner showed that they can be represented abstractly by certain recursively defined 
permutations called pats, and deduced that they are counted by the 
Catalan numbers. Counting pats by number of descents yields 
the curious identity
\[
\sum_{k=0}^{n}\frac{1}{2n-2k+1}\binom{2n-2k+1}{k}\binom{2k}{n-k} = 
C_{n},
\]
where only the middle third of the summands are nonzero.
\end{abstract}
\section{Introduction}
\vspace*{-5mm}
Martin Gardner \cite{gard} showed how to construct 
\htmladdnormallink{hexaflexagons}{http://en.wikipedia.org/wiki/Flexagon}
 in his 
1956 debut column in Scientific American. Soon after, a 
mathematical treatment by Oakley and Wisner appeared in the Monthly 
\cite{oakley}. They identified hexaflexagons with certain integer permutations
that they called pats (I don't know why). Pats are defined 
recursively, with permutations represented as lists. A singleton permutation 
is a pat, and a permutation $p$ of length $n\ge 2$ is a pat iff
(i) there is a unique split point that divides $p$ into 
subpermutations $p_{1},\ p_{2}$ such that all entries in $p_{1}$ are 
greater than all entries in $p_{2}$, and (ii) the reverse of each of 
$p_{1}$ and $p_{2}$ is a pat. The pats on $[3]$ are 231 and 312, and 
on [4] are 2431, 3241, 3412, 4132, 4213. The number of pats on 
$[n+1]$ (as Oakley and Wisner showed and will be evident) is the 
Catalan number $C_{n}=\frac{1}{n+1}\binom{2n}{n}$.

The purpose of this note is to find the distribution of the 
statistic ``number of descents'' on pats. (A descent in a permutation 
$p$ is a pair of adjacent entries $(p_{i},p_{i+1})$ with $p_{i}>p_{i+1}$.)
To do so, we establish a recurrence that leads to  a quartic equation
for the generating function. Lagrange inversion
then yields a simple closed form for the number of pats on 
$[n]$ with $k$ descents, thereby giving a combinatorial 
interpretation of the identity
\begin{equation}\label{main}
\sum_{k=0}^{n}\frac{1}{2n-2k+1}\binom{2n-2k+1}{k}\binom{2k}{n-k} = C_{n}.
\end{equation}
The last section exhibits a bijection from pats to full binary trees 
showing that descents on 
pats are distributed as even-level vertices on binary trees and 
concludes with 
a generalization of (\ref{main}). 
\section{Distribution of Descents on Pats}\label{distribution}

\vspace*{-5mm}
From its definition, a pat $p$ of length $n$ determines two pats 
$p_{1},p_{2}$ of lengths $i$ and $n-i$ respectively for some $i\in 
[1,n-1]$. The connecting entries from $p_{1}$ to $p_{2}$ contribute a 
descent to $p$ and the remaining descents of $p$ correspond to the 
ascents in $p_{1}$ and $p_{2}$. Hence we obtain the recurrence
\[
\hspace*{30mm}u(n,k)=\sum_{i=1}^{n-1}\sum_{j=0}^{i-1}u(i,j)u(n-i,n-j-k-1) 
\hspace*{20mm} n\ge 2,\ k\ge 0
\]
for the number $u(n,k)$ of pats on $[n]$ with $k$ descents, with 
$u(1,0)=1$. The recurrence leads directly to a functional equation for the generating function 
$F(x,y):=\sum_{n\ge 1,\,k\ge 0}u(n,k)x^{n}y^{k}$:
\begin{equation}\label{recur}
F(x,y)=x+\frac{\raisebox{-.5ex}{1}}{\raisebox{.5ex}{$y$}}\:F\big(xy,
\frac{\raisebox{-.5ex}{1}}{\raisebox{.5ex}{$y$}}\big)^{2}.
\end{equation}
Iterating (\ref{recur}) once yields the algebraic equation
\[
F=x+y\big(x+F^{2}\big)^{2},
\]
where $F(x,y)$ is now abbreviated to $F$. 
Direct application of  Lagrange inversion to solve this equation is 
cumbersome but a trick shown to me by Ira
Gessel (see \cite{multivariable}\,) rapidly 
solves it. Introduce a new variable $z$ and consider 
$F=F(x,y,z)$ defined by
\begin{equation}\label{ges}
F=z\Big(x+y\big(x+F^{2}\big)^{2}\Big)    
\end{equation}
Equation (\ref{ges}) has 
the form $z=F/\phi(F) $ where 
\[
\phi(F)=x+y\big(x+F^{2}\big)^{2},
\]
and so Lagrange inversion says that 
\begin{equation}\label{sol}
[\,y^{k}z^{2j+1}]F(x,y,z) 
=\frac{1}{\raisebox{.3ex}{\footnotesize{$2j+1$}}}\,[\,y^{k}F^{2j}]\phi(F)^{2j+1}=
\frac{1}{2j+1}\binom{2j+1}{k}\binom{2k}{j}x^{j+k+1}.
\end{equation}
The coefficient of $x^{n}y^{k}$ in $F(x,y,1)$ is obtained by setting 
$j=n-k-1$ in (\ref{sol}), yielding
\[
u(n,k)=[x^{n}y^{k}]F(x,y)=[x^{n}y^{k}]F(x,y,1)=\frac{1}{2n-2k-1}\binom{2n-2k-1}{k}\binom{2k}{n-k-1},
\]
and hence identity (\ref{main}) after replacing $n$ by $n+1$.
Here are the first few values of $u(n,k)$.
\[
\begin{array}{r|ccccccc}
	n^{\textstyle{\,\backslash \,k}} & 0 & 1 & 2 & 3 & 4 & 5 & 6 \\
\hline 
	1&   1 & 0 & 0 & 0 & 0 & 0 & 0 \\ 
 	2&    0 & 1 & 0 & 0 & 0 & 0 & 0 \\ 
	3&    0 & 2 & 0 & 0 & 0 & 0 & 0 \\ 
	4&   0 & 1 & 4 & 0 & 0 & 0 & 0 \\ 
	5&   0 & 0 & 12 & 2 & 0 & 0 & 0 \\ 
	6&   0 & 0 & 12 & 30 & 0 & 0 & 0 \\ 
	7&   0 & 0 & 4 & 100 & 28 & 0 & 0 \\
        8&    0 & 0 & 0 & 140 & 280 & 9 & 0 \\ 
        9&   0 & 0 & 0 & 90 & 980 & 360 & 0 \\ 
 	10&   0 & 0 & 0 & 22 & 1680 & 2940 & 220
\end{array}
\]
\centerline{ \ Table of values of $u(n,k)$, the distribution of descents on pats}
\section{Pats as Trees}
\vspace*{-5mm}
A pat on $[n+1]$ can be represented, using its successive split 
points, by a vertex-labeled full binary tree on $2n$ edges as 
illustrated below. 
\begin{center} 
\begin{pspicture}(-6,-1.9)(6,5.5)
\psset{unit=1cm} 
\rput(0,-1.2){\textrm{{\small full binary tree for pat $4\,7\,9\,8\,5\,6\,2\,1\,3,$}}}
\rput(0,-1.7){\textrm{{\small interior vertices at even level are enlarged}}}
\rput(-5,2.3){\textrm{{\footnotesize 4}}}
\rput(-3.5,4.3){\textrm{{\footnotesize 7}}}
\rput(-2,5.3){\textrm{{\footnotesize 9}}}
\rput(-1,5.3){\textrm{{\footnotesize 8}}}
\rput(0,4.3){\textrm{{\footnotesize 5}}}
\rput(1,4.3){\textrm{{\footnotesize 6}}}

\rput(2,3.3){\textrm{{\footnotesize 2}}}
\rput(3,3.3){\textrm{{\footnotesize 1}}} 
\rput(4.5,2.3){\textrm{{\footnotesize 3}}} 

\psdots(-5,2)(-3,1)(-1,2)(-2.5,3)(.5,3)(-3.5,4)(-1.5,4)(-2,5)(-1,5)
(0,4)(1,4)(0,0) (3.5,1)(2.5,2)(2,3)(3,3)(4.5,2)

\psdots[dotscale=2](0,0)(-1,2)(2.5,2)(-1.5,4)

\psline(-2,5)(-1.5,4)(-1,5) 
\psline(-3.5,4)(-2.5,3)(-1.5,4)
\psline(-5,2)(-3,1)(-1,2)
\psline(-2.5,3)(-1,2)
\psline(0,4)(.5,3)(1,4)
\psline(.5,3)(-1,2)(-3,1)(0,0)(3.5,1)

\psline(2,3)(2.5,2)(3,3)
\psline(2.5,2)(3.5,1)(4.5,2)

\gray{
\rput(0,-.3){\textrm{{\footnotesize 479856213}}}
\rput(-3.7,0.9){\textrm{{\footnotesize 479856}}}
\rput(-1.7,2){\textrm{{\footnotesize 79856}}}
\rput(-3,3){\textrm{{\footnotesize 798}}}
\rput(-1.9,4){\textrm{{\footnotesize 98}}}
\rput(.1,3.1){\textrm{{\footnotesize 56}}}

\rput(3.9,0.9){\textrm{{\footnotesize 213}}}
\rput(2.9,2){\textrm{{\footnotesize 21}}}
}
  
\end{pspicture}
\end{center}

In fact, the labels are unnecessary; they can be uniquely recovered 
from the underlying tree. So we have a bijection from pats on $[n+1]$ 
to full binary trees on $2n$ edges. Under this bijection, a descent in 
the pat shows up as an interior vertex at even level in the tree (where 
rain water would collect between the two leaves corresponding to the 
descent).

Thus we have, pruning the leaf edges in a full binary tree to obtain a binary tree,

\noindent \textbf{Theorem}\ \ \emph{ Descents on pats are distributed as even-level 
vertices in binary trees.}

\noindent Similar considerations for ternary and higher order trees yield a generalization 
of (\ref{main}):
\begin{equation*}
\sum_{k=0}^{n}\frac{1}{rn-rk+1}\binom{rn-rk+1}{k}\binom{rk}{n-k} = 
\frac{1}{rn+1}\binom{rn+1}{n}
\end{equation*}
for $r\ge 2$.

\vspace*{5mm}

\noindent \textbf{Acknowledgment} \ I thank Ira Gessel for suggesting 
the Lagrange inversion technique used in Section \ref{distribution}.

\end{document}